\documentclass [reqno]{amsart}
\usepackage {latexsym,epsfig}
\usepackage {amsmath,amssymb}
\usepackage [all]{xy}
\usepackage [latin1]{inputenc}

\newif\ifpdf
        \ifx\pdfoutput\undefined
        \pdffalse 
        \else
        \pdfoutput=1 
        \pdftrue
        \fi

\newcommand {\CC}{{\mathbb C}}
\newcommand {\NN}{{\mathbb N}}
\newcommand {\PP}{{\mathbb P}}
\newcommand {\QQ}{{\mathbb Q}}
\newcommand {\RR}{{\mathbb R}}
\newcommand {\ZZ}{{\mathbb Z}}

\newcommand {\calT}{{\mathcal T}}

\newcommand {\frM}{{\mathfrak M}}

\newcommand {\vdim}{\mathop {\rm vdim} \nolimits}
\newcommand {\res}{\mathop {\rm res} \nolimits}
\newcommand {\pt}{{\rm pt}}

\newcommand {\llangle}{\langle\!\langle}
\newcommand {\rrangle}{\rangle\!\rangle}
\newcommand {\vfc}[1]{[#1]^{\rm virt}}
\newcommand {\xydiag}[1]{\[ \xymatrix @M=5pt {#1} \]}
\newtheorem {theorem}{Theorem}[section]
\newtheorem {lemma}[theorem]{Lemma}
\newtheorem {proposition}[theorem]{Proposition}
\newtheorem {corollary}[theorem]{Corollary}

\newtheorem {convention}[theorem]{Convention}

\theoremstyle {definition}
\newtheorem {example}[theorem]{Example}
\theoremstyle {remark}
\newtheorem {remark}[theorem]{Remark}

\begin {document}


\title [Topological recursion relations in higher genus]{%
  Topological recursion relations and Gromov-Witten invariants in higher genus}
\author {Andreas Gathmann}
\address {Universit\"at Kaiserslautern, Fachbereich Mathematik,
          Postfach 3049, 67653 Kaiserslautern, Germany}
\email {gathmann@mathematik.uni-kl.de}

\begin {abstract}
  We state and prove a topological recursion relation that expresses any
  genus-$g$ Gromov-Witten invariant of a projective manifold with at least a
  $(3g-1)$-st power of a cotangent line class in terms of invariants with fewer
  cotangent line classes. For projective spaces, we prove that these relations
  together with the Virasoro conditions are sufficient to calculate the full
  Gromov-Witten potential. This gives the first computationally feasible way
  to determine the higher genus Gromov-Witten invariants of projective spaces.
\end {abstract}

\maketitle


Consider a moduli space $ \bar M_{g,n} (X,\beta) $ of $n$-pointed genus-$g$
stable maps of class $ \beta $ to a complex projective manifold $X$. For any $
1 \le i \le n $ let $ \psi_i \in A^1 (\bar M_{g,n} (X,\beta)) $ be the $i$-th
cotangent line class, i.e.\ the first Chern class of the line bundle on $ \bar
M_{g,n} (X,\beta) $ whose fibers are the cotangent spaces of the underlying
curves at the $i$-th marked point. Equations in the Chow group of $ \bar
M_{g,n} (X,\beta) $ that express products of cotangent line classes in terms of
boundary classes (i.e.\ classes on moduli spaces of \emph {reducible} stable
maps) are called topological recursion relations.

It has been proven recently by E. Ionel that any product of at least $g$
cotangent line classes on $ \bar M_{g,n} (X,\beta) $ is a sum of boundary
cycles \cite {I}. Unfortunately, the corresponding topological recursion
relations are not yet known explicitly for general $g$. The $ g \le 2 $ cases
can be found in \cite {Ge}. In theory, it should be possible to derive the
equations for other (at least low) values of $g$ from Ionel's work. As $g$
grows however, the terms in the topological recursion relations become very
complicated, and their number seems to grow exponentially. Consequently,
Ionel's result is barely useful for actual computations, although it is of
course very interesting from a theoretical point of view.

In this paper, we will prove a seemingly much weaker topological recursion
relation that expresses only a product of at least $ 3g-1 $ cotangent line
classes at the same point in terms of boundary cycles. The idea to obtain this
relation is simple: we just pull back the obvious relation $ \psi_1^{3g-1}=0 $
on $ \bar M_{g,1} $ to $ \bar M_{g,n} (X,\beta) $ along the forgetful map, and
keep track of the various pull-back correction terms in a clever way.

The result is a topological recursion relation that is extremely easy to
state and apply. To be precise, denote by $ \llangle \tau_{m_1} (\gamma_1)
\cdots \tau_{m_n} (\gamma_n) \rrangle_g $ the genus-$g$ Gromov-Witten
correlation function (i.e.\ the generating function for all invariants
containing $ ev_1^* \gamma_1 \cdot \psi_1^{m_1} \cdots ev_n^* \gamma_n \cdot
\psi_n^{m_n} $; see section \ref {trr} for details). Let $ \{ T_a \} $ be a
basis of the cohomology of $X$, and denote by $ \{ T^a \} $ the Poincar\'e-dual
basis. Then for any $ m \ge 0 $
  \[ \llangle \tau_{3g-1+m}(\gamma) \rrangle_g =
       \sum_{i+j=3g-2} \llangle \tau_m(\gamma) T_a \rrangle^i
                       \llangle \tau_j(T^a) \rrangle_g, \]
where we use the Einstein summation convention over $a$, and where the
auxiliary correlation functions $ \llangle \cdots \rrangle^i $ are defined
recursively by
  \[ \llangle \tau_m(\gamma_1) \gamma_2 \rrangle^i =
       \llangle \tau_{m+1}(\gamma_1) \gamma_2 \rrangle^{i-1}
       - \llangle \tau_m(\gamma_1) T_a \rrangle_0
         \llangle T^a \gamma_2 \rrangle^{i-1} \]
with the initial condition
  \[ \llangle \cdots \rrangle^0 = \llangle \cdots \rrangle_0. \]
Unlike other topological recursion relations, our relations involve neither
sums over graphs nor invariants of genus other than $g$ and 0. Moreover, the
auxiliary functions $ \llangle \cdots \rrangle^i $ are ``universal'' in the
sense that they do not depend on $g$. All this makes our relations very easy
and fast to apply.

The application that we have in mind in this paper is the Virasoro conditions
for the Gromov-Witten invariants of projective spaces. It has been proven
recently by Givental that the Gromov-Witten potential of a projective space
$ \PP^r $ satisfies an infinite series of differential equations, called the
Virasoro conditions \cite {Gi}. It is easily checked that these equations allow
for recursion over the genus and the number of marked points in the following
sense: given $ g>0 $, $ n \ge 1 $, cohomology classes $ \gamma_2,\dots,\gamma_n
\in A^*(X) $, and non-negative integers $ m_2,\dots,m_n $, the Virasoro
conditions can express linear combinations of invariants
  \[ \langle \tau_m (\gamma) \tau_{m_2} (\gamma_2) \cdots
       \tau_{m_n} (\gamma_n) \rangle_{g,d} \]
(where $ m \ge 0 $, $ \gamma \in A^*(X) $, and the degree $ d \ge 0 $ vary) in
terms of other invariants with either smaller genus, or the same genus and
smaller number of marked points. There is one such invariant for every choice
of $m$, i.e.\ $ r+1 $ invariants for every choice of $d$. There is however only
one non-trivial Virasoro condition for every $d$. Consequently, the Virasoro
conditions alone are not sufficient to compute the Gromov-Witten invariants.

This is where the topological recursion relations come to our rescue. By
inserting them into the Virasoro conditions, we can effectively bound the value
of $m$ in the set of unknown invariants above, leaving only the invariants with
$ 0 \le m < 3g-1 $. This way we arrive at infinitely many linear Virasoro
conditions (one for every choice of $d$) for only $ 3g-1 $ invariants. It is
now of course strongly expected that this system should be solvable, i.e.\ that
the coefficient matrix of this system of linear equations has maximal rank $
3g-1 $. We will show that this is indeed always the case. In fact, we will show
that \emph {any} choice of $ 3g-1 $ distinct non-trivial Virasoro conditions
leads to a system of linear equations that determines the invariants uniquely.
We do this by computing the determinant of the corresponding coefficient
matrix: if we pick the Virasoro conditions associated to the degrees $
d_0,\dots,d_{3g-2} $ and reduce the cotangent line powers by our topological
recursion relations, we arrive at a system of $ 3g-1 $ linear equations for $
3g-1 $ invariants whose determinant is simply
  \[ \frac {\prod_{i>j} (d_i-d_j)}{\prod_{i=1}^{3g-2} i!} \cdot
       \prod_{i=1}^{3g-2} \left( i+\frac 12 \right)^{3g-1-i}, \]
which is obviously always non-zero. Therefore the Virasoro recursion works,
i.e.\ we have found a constructive (and not too complicated) way to compute the
Gromov-Witten invariants of $ \PP^r $ in any genus.

One should note that the result for the determinant above is remarkably simple,
given the complicated structure of the Virasoro equations and our somewhat
arbitrary choice of topological recursion relations. It would be interesting to
see whether there is some deeper connection between the Virasoro conditions and
the topological recursion relations that explains this easy result. It would
also be interesting to extend our result to other Fano varieties.

We should also mention that our results have already been conjectured some time
ago by Eguchi and Xiong \cite {EX}. In their paper, they motivate our
topological recursion relations by arguments from string theory. Assuming that
these relations hold, Eguchi and Xiong use them together with the Virasoro
conditions to compute a few examples of higher genus Gromov-Witten invariants
of projective spaces. From this point of view one can regard our paper as
providing a solid mathematical footing for \cite {EX}.

The paper is organized as follows. In section \ref {trr} we will establish the
topological recursion relations mentioned above. We will then describe in
section \ref {virasoro} how to apply these results to the Virasoro conditions
to get systems of linear equations for the Gromov-Witten invariants of
projective spaces. The proof that these systems of equations are always
solvable (i.e.\ the computation of the determinant mentioned above) is given in
section \ref {determinant}. Finally, we will list some numbers obtained with
our method in section \ref {numbers}.

A C++ program that implements the algorithm of our paper and computes the
Gromov-Witten invariants of projective spaces can be obtained from the author
on request.


\section {Topological recursion relations} \label {trr}

The goal of this section is to prove the topological recursion relation stated
in the introduction. To do so, we will first compare certain cycles in the
moduli spaces of stable and prestable curves.

For any $ g,n \ge 0 $ let $ \bar\frM_{g,n} $ be the moduli space of complex
$n$-pointed genus-$g$ prestable curves, i.e.\ the moduli space of tuples $
(C,x_1,\dots,x_n) $ where $C$ is a nodal curve of arithmetic genus $g$, and the
$ x_i $ are distinct smooth points of $C$. This is a proper (but not separated)
smooth Artin stack of dimension $ 3g-3+n $. The open substack corresponding to
irreducible curves is denoted $ \frM_{g,n} $.

Recall that a prestable curve $ (C,x_1,\dots,x_n) $ is called stable if every
rational component has at least 3 and every elliptic component at least 1
special point, where the special points are the nodes and the marked points.
The open substack of $ \bar \frM_{g,n} $ corresponding to stable curves is
denoted $ \bar M_{g,n} $. It is a proper, separated, smooth Deligne-Mumford
stack. If $ 2g+n \ge 3 $ there is a stabilization morphism $ s: \bar \frM_{g,n}
\to \bar M_{g,n} $ that contracts every unstable component.

On any of the above moduli spaces and for any of the marked points $ x_i $ we
define the cotangent line class, denoted $ \psi_{x_i} $, to be the first Chern
class of the line bundle on the moduli space whose fiber at the point $
(C,x_1,\dots,x_n) $ is the cotangent space $ T^\vee_{C,x_i} $. Cotangent line
classes do not remain unchanged under stabilization --- they receive correction
terms from the locus of reducible curves where the marked point is on an
unstable component. Our first task is therefore to compare the cycles $ s^*
\psi_{x_i} $ and $ \psi_{x_i} $ on $ \bar \frM_{g,n} $ for $ g>0 $. For
simplicity, we will do this here only in the case $ n=1 $. As there is then
only one cotangent line class, we will simply write it as $ \psi $ (with no
index).

Let us define the ``correction terms'' that we will pick up when pulling back
cotangent line classes. For fixed $ g>0 $ and any $ k \ge 0 $ let $ M_k $ be
the product
  \[ M_k = \bar \frM_{g,1} \times \underbrace {
       \bar \frM_{0,2} \times \cdots \times \bar \frM_{0,2}
     }_{\mbox {$k$ copies}}. \]
Obviously, points of $ M_k $ correspond to a collection $
(C^{(0)},\dots,C^{(k)}) $ of prestable curves with $ C^{(0)} $ of genus $g$
and $ C^{(i)} $ of genus 0 for $ i>0 $, together with marked points $ x_0 $ on
$ C^{(0)} $ and $ x_i,y_i $ on $ C^{(i)} $ for $ i>0 $.

There is a natural proper gluing morphism $ \pi: M_k \to \bar \frM_{g,1} $ that
sends any point $ (C^{(0)},\dots,C^{(k)},x_0,x_1,y_1,\dots,x_k,y_k) $ to the
nodal 1-pointed curve $ (C^{(0)} \cup \cdots \cup C^{(k)},x_k) $, where $
C^{(i-1)} $ is glued to $ C^{(i)} $ for all $ i=1,\dots,k $ by identifying $
x_{i-1} $ with $ y_i $. The morphism $ \pi $ is an isomorphism onto its image
on the open subset $ \frM_{g,1} \times \frM_{0,2} \times \cdots \times
\frM_{0,2} $. Obviously, a generic point in the image of $ \pi $ is just a
curve with $ k+1 $ components aligned in a chain, with the first component
having genus $g$, and all others being rational. The marked point is always on
the last component.

For any collection $ \lambda=(\lambda_0,\dots,\lambda_k) $ of $ k+1 $
non-negative integers we denote by $ Z(\lambda) $ the cycle
  \[ \pi_* (\psi_{x_0}^{\lambda_0} \cdots \psi_{x_k}^{\lambda_k}
       \cdot [M_k]) \in A_*(\bar \frM_{g,1}) \]
in the Chow group of $ \bar \frM_{g,1} $ (for intersection theory on Artin
stacks we refer to \cite {K}). The cycle $ Z(\lambda) $ has pure codimension $
k+\lambda_0+\cdots+\lambda_k $ in $ \bar \frM_{g,1} $. It can be represented
graphically as \begin {center} \includegraphics {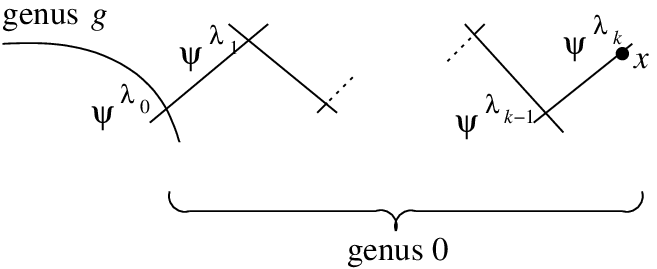} \end {center} where the
cotangent line classes sit at the ``left'' points of the nodes (except for the
last one that sits on the remaining marked point).

With these cycles $ Z(\lambda) $ we can now formulate the pull-back
transformation rule for cotangent line classes.

\begin {lemma} \label {spsi}
  For any $ k \ge 0 $ and $ \lambda_0,\dots,\lambda_k \ge 0 $ we have
    \[ s^* \psi \cdot Z(\lambda_0,\dots,\lambda_k) =
         Z(\lambda_0+1,\lambda_1,\dots,\lambda_k)
       - Z(0,\lambda_0,\lambda_1,\dots,\lambda_k) \]
  in $ A_* (\bar \frM_{g,1}) $, where $ s: \bar \frM_{g,1} \to \bar M_{g,1} $
  is the stabilization map.
\end {lemma}

\begin {proof}
  Let $ p_1: M_k = \bar \frM_{g,1} \times \bar \frM_{0,2} \times \cdots \times
  \bar \frM_{0,2} \to \bar \frM_{g,1} $ be the projection onto the first
  factor. Note that $ s \circ p_1 = s \circ \pi $.

  By \cite {Ge} proposition 5 we have the transformation rule
  \begin {center} \includegraphics {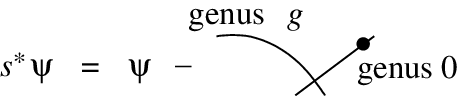} \end {center}
  on $ \bar \frM_{g,1} $. Pulling this equation back by $ p_1 $, intersecting
  with $ \psi_{x_0}^{\lambda_0} \cdots \psi_{x_k}^{\lambda_k} $, and pushing
  it forward again by $ \pi $ then gives the equation of the lemma.
\end {proof}

It is easy to iterate this lemma to get an expression for $ (s^* \psi)^N $:

\begin {corollary} \label {spsiN}
  For all $ N \ge 0 $ we have
    \[ (s^* \psi)^N = \sum_{k \ge 0} (-1)^k \sum_{\substack {
         (\lambda_0,\dots,\lambda_k) \\ k+\lambda_0+\cdots+\lambda_k = N
       }} Z(\lambda_0,\lambda_1,\dots,\lambda_k) \]
  in $ A_* (\bar \frM_{g,1}) $. In particular, the right hand side is zero if
  $ N \ge 3g-1 $.
\end {corollary}

\begin {proof}
  The statement is obvious for $ N=0 $ as $ \bar \frM_{g,1} = Z(0) $. The
  equation now follows immediately by induction from lemma \ref {spsi}.
  Moreover, note that $ \bar M_{g,1} $ is a Deligne-Mumford stack of dimension
  $ 3g-2 $, so its Chow groups in codimension at least $ 3g-1 $ vanish.
  Therefore $ (s^* \psi)^N = s^* (\psi^N) = 0 $ for $ N \ge 3g-1 $.
\end {proof}

\begin {remark} \label {spsiNpsim}
  Note that in the spaces $ M_k $ the marked point is always on the last
  component. So by intersecting the equation of corollary \ref {spsiN} with the
  $m$-th power of the cotangent line class on $ \bar \frM_{g,1} $ we get
    \[ \psi^m \cdot (s^* \psi)^N = \sum_{k \ge 0} (-1)^k \sum_{\substack {
         (\lambda_0,\dots,\lambda_k) \\ k+\lambda_0+\cdots+\lambda_k = N
       }} Z(\lambda_0,\lambda_1,\dots,\lambda_{k-1},\lambda_k+m) \]
  in $ A_* (\bar \frM_{g,1}) $ for all $ N,m \ge 0 $. As in the corollary, the
  right hand side will be zero if $ N \ge 3g-1 $.
\end {remark}

We will now apply this result to moduli spaces of stable maps. So let $X$ be a
complex projective manifold, and let $ \beta $ be the homology class of an
algebraic curve in $X$. As usual we denote by $ \bar M_{g,n} (X,\beta) $ the
moduli space of $n$-pointed genus-$g$ stable maps of class $ \beta $ to $X$
(see e.g.\ \cite {FP}). It is a proper Deligne-Mumford stack of virtual
dimension
  \[ \vdim \bar M_{g,n} (X,\beta) = -K_X\cdot \beta + (\dim X-3)(1-g) +n. \]
The actual dimension of $ \bar M_{g,n} (X,\beta) $ might (and for $ g>0 $
usually will) be bigger than this virtual dimension. There is however always a
canonically defined virtual fundamental class
  \[ \vfc {\bar M_{g,n} (X,\beta)} \in A_{\vdim \bar M_{g,n}(X,\beta)}
       (\bar M_{g,n} (X,\beta)) \]
that is used instead of the true fundamental class in intersection theory, and
that therefore makes the moduli space appear to have the ``correct'' dimension
for intersection-theoretic purposes (see e.g.\ \cite {BF}, \cite {B}).

The points in this moduli space can be written as $ (C,x_1,\dots,x_n,f) $,
where $ (C,x_1,\dots,x_n) \in \bar \frM_{g,n} $, $ f: C \to X $ is a morphism
of degree $ \beta $, and every rational (resp.\ elliptic) component \emph {on
which $f$ is constant} has at least 3 (resp.\ 1) special points.

The moduli spaces $ \bar M_{g,n} (X,\beta) $ come equipped with cotangent line
classes $ \psi_{x_i} $ in the same way as above. In addition, for every marked
point $ x_i $ we have an evaluation morphism $ ev_{x_i} : \bar M_{g,n}
(X,\beta) \to X $ given by $ ev_{x_i} (C,x_1,\dots,x_n,f) = f(x_i) $. For any
cohomology classes $ \gamma_1,\dots,\gamma_n \in A^*(X) $ and non-negative
integers $ m_1,\dots,m_n $ we define the Gromov-Witten invariant
  \[ \langle \tau_{m_1}(\gamma_1) \cdots \tau_{m_n} (\gamma_n)
     \rangle_{g,\beta} :=
     \int_{\vfc {\bar M_{g,n}(X,\beta)}}
       \psi_{x_1}^{m_1} \cdot ev_{x_1}^* \gamma_1 \cdots
       \psi_{x_n}^{m_n} \cdot ev_{x_n}^* \gamma_n \in \QQ, \]
where the integral is understood to be zero if the integrand is not of
dimension $ \vdim \bar M_{g,n} (X,\beta) $. If $ m_i = 0 $ for some $i$ we
abbreviate $ \tau_{m_i} (\gamma_i) $ as $ \gamma_i $ within the brackets on the
left hand side. As the Gromov-Witten invariants are multilinear in the $
\gamma_i $, it suffices to pick the $ \gamma_i $ from among a fixed basis.
So let us choose a basis $ \{ T_a \} $ of the cohomology (modulo numerical
equivalence) of $X$ and let $ \{ T^a \} $ be the Poincar\'e-dual basis.

\begin {remark} \label {corr1}
  It is often convenient to encode the Gromov-Witten invariants as the
  coefficients of a generating function. So we introduce the so-called
  correlation functions
    \[ \llangle
         \tau_{m_1}(\gamma_1) \cdots \tau_{m_n}(\gamma_n)
       \rrangle_g := \sum_\beta \left\langle
         \tau_{m_1}(\gamma_1) \cdots \tau_{m_n}(\gamma_n)
         \exp \left(\sum_m t_m^a \tau_m (T_a)\right) \right\rangle_{g,\beta}
         \; q^\beta
    \]
  where the $ t_m^a $ and $ q^\beta $ are formal variables satisfying $
  q^{\beta_1} q^{\beta_2} = q^{\beta_1+\beta_2} $. Here and in the following we
  use the summation convention for the ``cohomology index'' $a$, i.e.\ an index
  occurring both as a lower and upper index is summed over. The correlation
  functions are formal power series in the variables $ t_m^a $ and $ q^\beta $
  whose coefficients describe all genus-$g$ Gromov-Witten invariants containing
  at least the classes $ \tau_{m_1}(\gamma_1) \cdots \tau_{m_n}(\gamma_n) $.
\end {remark}

With this notation we can now rephrase corollary \ref {spsiN} in terms of
correlation functions for Gromov-Witten invariants:

\begin {proposition} \label {trr-prep}
  For all $ g>0 $, $ N \ge 3g-1 $, $ m \ge 0 $, and $ \gamma \in A^*(X) $ we
  have
  \begin {align*}
    0 &= \sum_{k \ge 0} (-1)^k
      \sum_{\substack {(\lambda_0,\dots,\lambda_k) \\
        k+\lambda_0+\cdots+\lambda_k = N }} \\
    &\qquad \qquad
    \llangle \tau_{\lambda_0} (T^{a_1}) \rrangle_g
    \llangle T_{a_1} \tau_{\lambda_1} (T^{a_2}) \rrangle_0
    \cdots
    \llangle T_{a_{k-1}} \tau_{\lambda_{k-1}} (T^{a_k}) \rrangle_0
    \llangle T_{a_k} \tau_{\lambda_k+m} (\gamma) \rrangle_0
  \end {align*}
  as power series in $ t_m^a $ and $ q^\beta $.
\end {proposition}

\begin {proof}
  For every $ n \ge 1 $ and any homology class $ \beta $ there is a forgetful
  morphism
    \[ q: \bar M_{g,n} (X,\beta) \to \bar \frM_{g,1}, \quad
          (C,z_1,\dots,z_n,f) \mapsto (C,z_1). \]
  (We denote the marked points by $ z_i $ in order not to confuse them with the
  $ x_i $ and $ y_i $ above that are used to glue the components of the
  reducible curves.) We claim that the statement of the proposition is obtained
  by pulling back the equation of remark \ref {spsiNpsim} by $q$ in the case $
  N \ge 3g-1 $ and evaluating the result on the virtual fundamental class of
  $ \bar M_{g,n} (X,\beta) $.

  To show this, we obviously have to compute $ q^* Z(\lambda) $ for all $
  \lambda = (\lambda_0,\dots,\lambda_k) $. Let $ B = (\beta_0,\dots,
  \beta_k) $ be a collection of homology classes with $ \sum_i \beta_i =
  \beta $, and let $ I = (I_0,\dots,I_k) $ be a collection of subsets of $
  \{2,\dots,n\} $ whose union is $ \{2,\dots,n\} $. Set
    \[ M_{B,I} = \bar M_{g,1+\#I_0} (X,\beta_0) \times_X
                 \bar M_{0,2+\#I_1} (X,\beta_1) \times_X \cdots \times_X
                 \bar M_{0,2+\#I_k} (X,\beta_k), \]
  where the fiber products are taken over the evaluation maps at the first
  marked point of the $(i-1)$-st factor and the second marked point of the
  $i$-th factor for $ i=1,\dots,k $. In other words, the moduli space $ M_{B,I}
  $ describes stable maps with the same configuration of components as in $ M_k
  $, and with the homology class $ \beta $ and the marked points $
  z_2,\dots,z_n $ split up onto the components in a prescribed way (the point $
  z_1 $ is always the first marked point of the last factor). Note that the
  space $ M_{B,I} $ carries a natural virtual fundamental class induced from
  their factors. By \cite {B} axiom III it is equal to the product cycle of the
  virtual fundamental classes of the factors, intersected with the pull-backs
  of the diagonal classes $ \Delta_X \subset X \times X $ along the evaluation
  maps at every pair of marked points where two components are glued together.

  Now by \cite {B} axiom V we have a Cartesian diagram
  \xydiag {
    \coprod_{B,I} M_{B,I} \ar[r]^{\tilde \pi} \ar[d]_{\tilde q} &
      \bar M_{g,n} (X,\beta) \ar[d]^q \\
    M_k \ar[r]^{\pi} &
      \frM_{g,1}
  }
  with $ \pi^! \vfc {\bar M_{g,n}(X,\beta)} = \sum_{B,I} \vfc {M_{B,I}} $. As
  $ \tilde q $ does not change the cotangent line classes we conclude that
  \begin {align*}
    q^* Z(\lambda) \cdot \vfc {\bar M_{g,n} (X,\beta)}
    &= q^* \pi_* (\psi_{x_0}^{\lambda_0} \cdots
        \psi_{x_k}^{\lambda_k}) \cdot \vfc {\bar M_{g,n} (X,\beta)} \\
    &= {\tilde \pi}_* {\tilde q}^* (\psi_{x_0}^{\lambda_0} \cdots
        \psi_{x_k}^{\lambda_k}) \cdot \vfc {\bar M_{g,n} (X,\beta)} \\
    &= \sum_{B,I} {\tilde \pi}_* \left(
         \psi_{x_0}^{\lambda_0} \cdots \psi_{x_k}^{\lambda_k}
	 \cdot \vfc {M_{B,I}} \right), \tag {$*$}
  \end {align*}
  where $ x_i $ denotes the first marked point of the $i$-th factor. Now
  choose cohomology classes $ \gamma_1,\dots,\gamma_n $ and non-negative
  integers $ m_1,\dots,m_n $. Intersecting expression $ (*) $ with $
  ev_{z_1}^* \gamma_1 \cdot \psi_{z_1}^{m_1} \cdots ev_{z_n}^* \gamma_n \cdot
  \psi_{z_n}^{m_n} $ and taking the degree of the resulting homology class (if
  it is zero-dimensional), we get exactly the Gromov-Witten invariants
    \[ \sum_{B,I} \langle \tau_{\lambda_0} (T^{a_1})
         \calT_0 \rangle_{g,\beta_0}
       \langle T_{a_1} \tau_{\lambda_1} (T^{a_2})
         \calT_1 \rangle_{0,\beta_1}
       \cdots
       \langle T_{a_k} \tau_{\lambda_k+l_1} (\gamma_1)
         \calT_k \rangle_{0,\beta_k}, \]
  where $ \calT_i $ is a short-hand notation for $ \prod_{j \in I_i} \tau_{m_j}
  (\gamma_j) $. If we choose $ \gamma_i = T^{a_i} $ for some $ a_i $, this can
  obviously be rewritten as the $ \left( q^\beta \cdot \prod_{i=1}^n
  t_{m_i}^{a_i} \right) $-coefficient of the function
    \[ \llangle \tau_{\lambda_0} (T^{a_1}) \rrangle_g
       \llangle T_{a_1} \tau_{\lambda_1} (T^{a_2}) \rrangle_0
       \cdots
       \llangle T_{a_{k-1}} \tau_{\lambda_{k-1}} (T^{a_k}) \rrangle_0
       \llangle T_{a_k} \tau_{\lambda_k+m} (\gamma) \rrangle_0. \]
  Inserting this into the formula of remark \ref {spsiNpsim} gives the desired
  result.
\end {proof}

\begin {corollary} \label {trr-final} \textbf {(Topological recursion
  relation)}
  For all $ g>0 $, $ N \ge 3g-1 $, $ m \ge 0 $, and $ \gamma \in A^*(X) $ we
  have
    \[ \llangle \tau_{N+m}(\gamma) \rrangle_g =
         \sum_{i+j=N-1} \llangle \tau_m(\gamma) T_a \rrangle^i
                         \llangle \tau_j(T^a) \rrangle_g, \]
  where the auxiliary correlation functions $ \llangle \cdots \rrangle^i $ are
  defined recursively by
    \[ \llangle \tau_m(\gamma_1) \gamma_2 \rrangle^i =
         \llangle \tau_{m+1}(\gamma_1) \gamma_2 \rrangle^{i-1}
         - \llangle \tau_m(\gamma_1) T_a \rrangle_0
           \llangle T^a \gamma_2 \rrangle^{i-1} \]
  with the initial condition
    \[ \llangle \cdots \rrangle^0 = \llangle \cdots \rrangle_0. \]
\end {corollary}

\begin {proof}
  Note that the $ k=0 $ term in lemma \ref {trr-prep} is just $ \llangle
  \tau_{N+m}(\gamma) \rrangle_g $. So we find that the equation of the
  corollary is true if we set
  \begin {align*}
    \llangle \tau_m(\gamma_1) \gamma_2 \rrangle^i
    &= \sum_{k \ge 0} (-1)^k
      \sum_{\substack {(\lambda_1,\dots,\lambda_k) \\
        k+\lambda_1+\cdots+\lambda_k = N }} \\
    &\qquad \qquad
    \llangle T_{a_1} \tau_{\lambda_1} (T^{a_2}) \rrangle_0
    \cdots
    \llangle T_{a_{k-1}} \tau_{\lambda_{k-1}} (T^{a_k}) \rrangle_0
    \llangle T_{a_k} \tau_{\lambda_k+m} (\gamma) \rrangle_0.
  \end {align*}
  It is checked immediately that these correlation functions satisfy the
  recursive relations stated in the corollary.
\end {proof}

\begin {remark} \label {corr2}
  As in the case of the Gromov-Witten invariants, we will expand the
  correlation functions $ \llangle \cdots \rrangle^i $ as a power series in $
  q^\beta $ and $ t_m^a $ and call the resulting coefficients $ \langle
  \cdots \rangle^i_\beta $ according to the formula
    \[ \llangle
         \tau_{m_1}(\gamma_1) \gamma_2
       \rrangle^i = \sum_\beta \left\langle
         \tau_{m_1}(\gamma_1) \gamma_2
         \exp \left(\sum_m t_m^a \tau_m (T_a)\right) \right\rangle^i_\beta
	 \; q^\beta.
    \]
  Note however that, in contrast to the Gromov-Witten numbers, the invariants
  $ \langle \cdots \rangle^i $ must have at least two entries, of which the
  second one contains no cotangent line class.
\end {remark}

For future computations it is convenient to construct a minor generalization of
corollary \ref {trr-final} that is mostly notational. Note that all genus-0
degree-0 invariants with fewer than 3 marked points are trivially zero, as
the moduli spaces of stable maps are empty in this case. It is an important and
interesting fact that many formulas concerning Gromov-Witten invariants get
easier if we assign ``virtual values'' to these invariants in the unstable
range:

\begin {convention} \label {conv}
  Unless stated otherwise, we will from now on allow formal negative powers of
  the cotangent line classes (i.e.\ the index $m$ in the $ \tau_m(\gamma) $ can
  be any integer). Invariants $ \langle \cdots \rangle_{g,\beta} $ and $
  \langle \cdots \rangle^i_\beta $ are simply defined to be zero if they
  contain a negative power of a cotangent line class at any point, except for
  the following cases of genus-0 degree-0 invariants with fewer than 3 marked
  points:
  \begin {enumerate}
  \item $ \langle \tau_{-2} (\pt) \rangle_{0,0} = 1 $,
  \item $ \langle \tau_{m_1} (\gamma_1) \tau_{m_2} (\gamma_2) \rangle_{0,0} =
    (-1)^{\max (m_1,m_2)}(\gamma_1 \cdot \gamma_2) \delta_{m_1+m_2,-1} $,
  \item $ \langle \tau_{-i-1} (\gamma_1) \gamma_2 \rangle^i_0 = (\gamma_1 \cdot
    \gamma_2) $ for all $ i \ge 0 $.
  \end {enumerate}
  The correlation functions $ \llangle \cdots \rrangle $ are changed
  accordingly so that the equations of remarks \ref {corr1} and \ref {corr2}
  remain true (in particular these functions will now depend additionally on
  the variables $ t_m^a $ for $ m<0 $).
\end {convention}

\begin {remark}
  Note that this convention is consistent with the general formula for genus-0
  degree-0 invariants
    \[ \langle \tau_{m_1}(\gamma_1) \cdots \tau_{m_n} (\gamma_n)
       \rangle_{0,0} =
         \binom {n-3}{m_1,\dots,m_n} (\gamma_1 \cdots \gamma_n)
	 \delta_{m_1+\cdots+m_n,n-3}, \]
  as well as with the recursion relations for the $ \llangle \cdots \rrangle^i
  $ of corollary \ref {trr-final}.
\end {remark}

Using this convention, we can now restate our topological recursion relations
as follows:

\begin {corollary} \label {trr-conv}
  For all $ g>0 $, $ N \ge 3g-1 $, $ m \in \ZZ $, and $ \gamma \in A^*(X) $ we
  have
    \[ \llangle \tau_{N+m} (\gamma) \rrangle_g =
         \sum_{i+j=N-1} \llangle \tau_m(\gamma) T_a \rrangle^i
                         \llangle \tau_j(T^a) \rrangle_g, \]
  where the auxiliary correlation functions $ \llangle \cdots \rrangle^i $ are
  defined recursively by the formulas given in corollary \ref {trr-final},
  together with convention \ref {conv}.
\end {corollary}

\begin {proof}
  The equations in the corollary are the same as in corollary \ref {trr-final}
  if $ m \ge 0 $. For $ m<0 $ they reduce to the trivial equations $ \llangle
  \tau_{N+m} (\gamma) \rrangle_g = \llangle \tau_{N+m} (\gamma) \rrangle_g $ by
  convention \ref {conv}.
\end {proof}


\section {The Virasoro conditions} \label {virasoro}

We now want to apply our topological recursion relation in conjunction with the
Virasoro conditions to compute Gromov-Witten invariants. The Virasoro
conditions are certain relations among Gromov-Witten invariants conjectured in
\cite {EHX} that have recently been proven for projective spaces by Givental
\cite {Gi}. We will therefore from now on restrict to the case $ X=\PP^r $. It
is expected that the same methods would work for other Fano varieties as well.

To state the Virasoro conditions we need some notation. We pick the obvious
basis $ \{T_a\} $ of $ A^*(X) $ where $ T_a $ denotes the class of a linear
subspace of codimension $a$ for $ a=0,\dots,r $. Let $ R: A^*(X) \to A^*(X) $
be the homomorphism of multiplication with the first Chern class $ c_1(X) $. In
our basis, the $p$-th power $ R^p $ of $R$ is then given by $ (R^p)_a{}^b =
(r+1)^p \, \delta_{a+p,b} $.

For any $ x \in \QQ $, $ k \in \ZZ_{\ge-1} $, and $ 0 \le p \le k+1 $ denote
by $ [x]^k_p $ the $ z^p $-coefficient of $ \prod_{j=0}^k (z+x+j) $, or in
other words the $ (k+1-p) $-th elementary symmetric polynomial in $ k+1 $
variables evaluated at the numbers $ x,\dots,x+k $.

Then the Virasoro conditions state that for any $ k \ge 1 $ and $ g \ge 1 $ we
have an equation of power series in $ t_m^a $ and $ q^\beta $ (see e.g.\
\cite {EHX})
\begin {small} \begin {align*}
  0 =& - \sum_{p=0}^{k+1} \left[ \frac {3-r}2 \right]^k_p (R^p)_0{}^b
         \llangle \tau_{k+1-p} (T_b) \rrangle_g \tag {A} \\
  & + \sum_{p=0}^{k+1} \sum_{m=0}^\infty
         \left[ a+m+\frac {1-r}{2} \right]^k_p (R^p)_a{}^b t_m^a
         \llangle \tau_{k+m-p} (T_b) \rrangle_g \tag {B} \\
  & + \frac 12 \sum_{p=0}^{k+1} \sum_{m=p-k}^{-1} (-1)^m
         \left[ a+m+\frac {1-r}{2} \right]^k_p (R^p)_a{}^b
         \llangle \tau_{-m-1}(T^a) \tau_{k+m-p}(T_b) \rrangle_{g-1} \tag {C} \\
  & + \frac 12 \sum_{p=0}^{k+1} \sum_{m=p-k}^{-1} \sum_{h=0}^g (-1)^m
         \left[ a+m+\frac {1-r}{2} \right]^k_p (R^p)_a{}^b
         \llangle \tau_{-m-1}(T^a) \rrangle_h
         \llangle \tau_{k+m-p}(T_b) \rrangle_{g-h}, \tag {D}
\end {align*} \end {small}%
where convention \ref {conv} is not yet applied (i.e.\ the genus-0 degree-0
invariants in the unstable range are defined to be zero). We should also
mention that there are versions of these relations also for $ k \ge -1 $ and
all $ g \ge 0 $, but the equations will then get additional correction terms
that we have dropped here for the sake of simplicity.

First of all let us apply convention \ref {conv} to these formulas. It is
checked immediately that this realizes the (A) and (B) terms as part of the (D)
terms via the conventions (i) and (ii), respectively. So by applying our
convention we can drop the (A) and (B) terms above (if we allow arbitrary
integers in the sum over $m$).

Let us now analyze how these equations can be used to compute Gromov-Witten
invariants. First of all we will compute the invariants recursively over the
genus of the curves. The genus-0 invariants of $ \PP^r $ are well-known and can
be computed by the WDVV equations (see e.g.\ the ``first reconstruction
theorem'' of Kontsevich and Manin \cite {KM}). So let us assume that we want to
compute the invariants of some genus $ g>0 $, and that we already know all
invariants of smaller genus. In the Virasoro equations above this means that we
know all of (C), as well as the terms of (D) where $ h \neq 0 $ and $ h \neq g
$. Noting that the (D) terms are symmetric under $ h \mapsto g-h $, we can
therefore rewrite the Virasoro conditions as
\begin {align*}
  & \sum_{p=0}^{k+1} \sum_m (-1)^m
      \left[ a+m+\frac {1-r}{2} \right]^k_p (R^p)_a{}^b
      \llangle \tau_{-m-1}(T^a) \rrangle_0
      \llangle \tau_{k+m-p}(T_b) \rrangle_g \\
  & \qquad \qquad = \mbox {(recursively known terms)}.
\end {align*}
Next, we will compute the invariants of genus $g$ recursively over the number
of marked points. So let us assume that we want to compute the $n$-point
genus-$g$ invariants, and that we already know all invariants of genus $g$ with
fewer marked points. In the equations above this means that we fix a degree $
d \ge 0 $, integers $ m_2,\dots,m_n $, and $ n-1 $ cohomology classes $
T_{a_2},\dots,T_{a_n} $, and compare the $ \left( q^d \cdot \prod_{i=2}^n
t^{a_i}_{m_i} \right) $-coefficients of the equations. By the recursion process
we then know all the invariants in which at least one of the marked points $
x_2,\dots,x_n $ is on the genus-0 invariant. So we can write
\begin {align*}
  & \sum_{p=0}^{k+1} \sum_m \sum_{d_1+d_2=d} (-1)^m
       \left[ a+m+\frac {1-r}{2} \right]^k_p (R^p)_a{}^b \cdot \\
  & \qquad \qquad \qquad \qquad \qquad \cdot
       \langle \tau_{-m-1}(T^a) \rangle_{0,d_1}
       \langle \tau_{k+m-p}(T_b)
       \tau_{m_2}(T_{a_2}) \cdots \tau_{m_n}(T_{a_n})
       \rangle_{g,d_2} \\
  & \qquad \qquad = \mbox {(recursively known terms)}.
\end {align*}
These are equations for the unknown invariants $ \langle \tau_j(T_b)
\tau_{m_2}(T_{a_2}) \cdots \tau_{m_n}(T_{a_n})\rangle_{g,e} $, where $j$, $b$,
and $e$ vary. Note that for a given $ j \ge 0 $ there is exactly one such
invariant $ \langle \tau_j(T_{b_j}) \tau_{m_2}(T_{a_2}) \cdots \tau_{m_n}
(T_{a_n}) \rangle_{g,e_j} $: the values of $ b_j $ and $ e_j $ are determined
uniquely by the dimension condition
\begin {equation} \label {dim-cond}
  (r+1)e_j + (r-3)(1-g)+n = j + b_j + \sum_{i=2}^n (m_i+a_i)
\end {equation}
as we must have $ 0 \le b_j \le r $. Let us denote this invariant by $ x_j $.
Of course it may happen that $ e_j<0 $, in which case we set $ x_j=0 $. Our
equations now read
\begin {align*}
  & \sum_{p=0}^{k+1} \sum_m (-1)^{m+p-k}
       \left[ a+m+p-k+\frac {1-r}{2} \right]^k_p (R^p)_a{}^{b_m}
       \langle \tau_{-m-p+k-1}(T^a) \rangle_{0,d-e_m} \cdot x_m \\
  & \qquad \qquad = \mbox {(recursively known terms)}.
\end {align*}
Let us now check how many non-trivial equations of this sort we get. Together
with (\ref {dim-cond}), the dimension conditions
  \[ (d-e_m)(r+1)+r-3+1 = -m-p+k-1+r-a \]
(for the genus-0 invariant) and $ a+p=b_m $ (from the $ R^p $ factor) give
\begin {equation} \label {k-det}
  k = d(r+1) + (r-3)(1-g) +n-1-\sum_{i=2}^n (m_i+a_i),
\end {equation}
which means that the value of $k$ is determined by $d$. To avoid overly
complicated notation, in what follows we will denote the number $k$ determined
by (\ref {k-det}) by $ k(d) $. Moreover, let $ \delta $ be the smallest value
of $d$ for which $ k(d) $ is positive. We are then getting one equation for
every degree $ d \ge \delta $. As there are $ r+1 $ unknown invariants $ x_j $
in every degree however, it is clear that our equations alone are not
sufficient to determine the $ x_j $.

Let us now apply our topological recursion relations. In terms of the recursion
at hand, these relations can express every invariant $ x_m $ as a linear
combination of invariants of the same form with $ m<3g-1 $, plus some terms
that are known recursively because they contain only invariants with fewer than
$n$ marked points. More precisely, we have
  \[ x_m = \sum_{i+j=N-1}
       \langle \tau_{m-N+2} (T_{b_m}) T^{b_j} \rangle^i_{e_m-e_j} x_j
       + \mbox {(recursively known terms)} \]
for all $ N \ge 3g-1 $ by corollary \ref {trr-conv}. Inserting this into the
Virasoro conditions, we get
\begin {align*}
  & \sum_{p=0}^{k(d)+1} \sum_m \sum_{i+j=N-1} (-1)^{m+p-k(d)}
       \left[ a+m+p-k(d)+\frac {1-r}{2} \right]^{k(d)}_p (R^p)_a{}^{b_m}
       \cdot \\
  & \qquad \qquad \qquad \qquad \cdot
       \langle \tau_{-m-p+k(d)-1}(T^a) \rangle_{0,d-e_m} \cdot
       \langle \tau_{m-N+2} (T_{b_m}) T^{b_j}
       \rangle^i_{e_m-e_j} \cdot
       x_j \\
  & \qquad \qquad = \mbox {(recursively known terms)}.
\end {align*}
Using the dimension conditions again, and noting that the sum over $m$ is
equivalent to independent sums over $ b_m $ and $ e_m $, we can rewrite this as
\begin {align*}
  & \sum_{p=0}^{k(d)+1} \sum_e \sum_{i+j=N-1} (-1)^{1-a-(d-e)(r+1)}
       \left[ \frac {3-r}{2}-(d-e)(r+1) \right]^{k(d)}_p (R^p)_a{}^b \cdot \\
  & \qquad \qquad \qquad \qquad \cdot
       \langle \tau_{\bullet}(T^a) \rangle_{0,d-e} \cdot
       \langle \tau_{\bullet} (T_b) T^{b_j}
       \rangle^i_{e-e_j} \cdot
       x_j \\
  & \qquad \qquad = \mbox {(recursively known terms)},
\end {align*}
where the dots in the $ \tau $ functions denote the uniquely determined numbers
so that the invariants satisfy the dimension condition.

We are thus left with infinitely many equations (one for every $ d \ge \delta
$) for finitely many variables $ x_0,\dots,x_{N-1} $. It is of course strongly
expected that this system of equations should be solvable, i.e.\ that the
matrix $ V^{(N)}=(V^{(N)}_{d,j})_{d \ge \delta,0 \le j < N} $ with
\begin {align} \label {v-def}
  & V^{(N)}_{d,j} := \sum_{p=0}^{k(d)+1} \sum_e (-1)^{1-a-(d-e)(r+1)}
       \left[ \frac {3-r}{2}-(d-e)(r+1) \right]^{k(d)}_p (R^p)_a{}^b \cdot
       \notag \\
  & \qquad \qquad \qquad \qquad \cdot
       \langle \tau_{\bullet}(T^a) \rangle_{0,d-e} \cdot
       \langle \tau_{\bullet}(T_b) T^{b_j}
       \rangle^{N-1-j}_{e-e_j}
\end {align}
has maximal rank $N$. This is what we will show in the next section. In fact,
we will prove that \emph {every} $ N \times N $ submatrix of $V$ is invertible.
So we have shown

\begin {theorem} \label {reconstruction}
  The Virasoro conditions together with the topological recursion relations of
  corollary \ref {trr-conv} give a constructive way to determine all
  Gromov-Witten invariants of projective spaces.
\end {theorem}

In contrast to other known relations that in theory determine the Gromov-Witten
invariants (see \cite {GP}, \cite {LLY}), our algorithm is easily implemented
on a computer. No complicated sums over graphs occur anywhere in the
procedure. It should be noted however that the calculation of some genus-$g$
degree-$d$ invariant usually requires the recursive calculation of invariants
of smaller genus with \emph {bigger} degree and \emph {more} marked points.
This is the main factor for slowing down the algorithm as the genus grows.

Some numbers that have been computed using this algorithm can be found in
section \ref {numbers}.


\section {Computation of the determinant} \label {determinant}

The goal of this section is to prove the technical result needed for theorem
\ref {reconstruction}:

\begin {proposition} \label {det-nonzero}
  Fix any $ N \ge 1 $, and let $ V^{(N)}=(V^{(N)}_{d,j})_{d \ge \delta,0 \le
  j<N} $ be the matrix defined in equation (\ref {v-def}). Then any $ N \times
  N $ submatrix of $ V^{(N)} $, obtained by picking $N$ distinct values of $d$,
  has non-zero determinant.
\end {proposition}

We will prove this statement in several steps. In a first step, we will make
the entries of the matrix independent of $N$ and reduce the invariants $
\langle \cdots \rangle^i $ to ordinary rational Gromov-Witten invariants:

\begin {lemma}
  Let $ W=(W_{d,j})_{d \ge \delta,j \ge 0} $ be the matrix with entries $
  W_{d,j} = V^{(j+1)}_{d,j} $. Then:
  \begin {enumerate}
  \item For all $ d \ge \delta $, $ N \ge 1 $, and $ 0 \le j < N $ we have
      \[ V_{d,j}^{(N+1)} = V_{d,j}^{(N)}
           - \langle T_{b_N} T^{b_j} \rangle^{N-1-j}_{e_N-e_j} \cdot W_{d,N}.
      \]
  \item For any $ N \ge 1 $ and any $ N \times N $ submatrix of $W$ obtained
    by taking the \emph {first} $N$ columns of \emph {any} $N$ rows, the
    determinant of this submatrix is the same as the corresponding submatrix
    of $ V^{(N)} $.
  \end {enumerate}
\end {lemma}

\begin {proof}
  (i): Comparing the $ q^{e-e_j} $-terms of the recursive relations of
  corollary \ref {trr-final} we find that
    \[ \langle \tau_{\bullet}(T_b) T^{b_j} \rangle^{N-j}_{e-e_j}
       = \langle \tau_{\bullet}(T_b) T^{b_j} \rangle^{N-1-j}_{e-e_j}
         - \langle \tau_{\bullet}(T_b) T^{b_N} \rangle_{0,e-e_N}
           \langle T_{b_N} T^{b_j} \rangle^{N-1-j}_{e_N-e_j}, \]
  from which the claim follows.

  (ii): We prove the statement by induction on $N$. There is nothing to show
  for $ N=1 $. Now assume that we know the statement for some value of $N$,
  i.e.\ any two corresponding $ N \times N $ submatrices of the matrices with
  columns
    \[ (W_{\cdot,0},\dots,W_{\cdot,N-1}) \quad \mbox {and} \quad
       (V^{(N)}_{\cdot,0},\dots,V^{(N)}_{\cdot,N-1}) \]
  have the same determinant. Of course, the same is then also true for any
  corresponding $ (N+1) \times (N+1) $ submatrices of
    \[ (W_{\cdot,0},\dots,W_{\cdot,N-1},W_{\cdot,N}) \quad \mbox {and} \quad
       (V^{(N)}_{\cdot,0},\dots,V^{(N)}_{\cdot,N-1},W_{\cdot,N}). \]
  But by (i), the latter matrix is obtained from
    \[ (V^{(N+1)}_{\cdot,0},\dots,V^{(N+1)}_{\cdot,N-1},W_{\cdot,N}) =
       (V^{(N+1)}_{\cdot,0},\dots,V^{(N+1)}_{\cdot,N}) \]
  by an elementary column operation, so the result follows.
\end {proof}

So by the lemma, it suffices to consider the matrix $W$. Let us now evaluate
the genus-0 Gromov-Witten invariants contained in the definition of $W$.

\begin {convention} \label {prod-conv}
  For the rest of this section, we will make the usual convention that a
  product $ \prod_{i=i_1}^{i_2} A_i $ is defined to be $
  \prod_{i=i_2+1}^{i_1-1} A_i^{-1} $ if $ i_1>i_2 $.
\end {convention}

\begin {lemma}
  For all $ d \ge \delta $ and $ j \ge 0 $ the matrix entry $ W_{d,j} $ is
  equal to the $ z^j $-coefficient of
    \[ - \frac {\prod_{i=0}^{k(d)} \left(
         r+1+\left( \frac {3-r}2+i\right) z \right)}{
         (1+(d-e_j)z)^{b_j+1} \prod_{i=0}^{d-e_j-1} (1+iz)^{r+1}}. \]
\end {lemma}

\begin {proof}
  Recall that by equation (\ref {v-def}) the matrix entries $ W_{d,j} $ are
  given by
  \begin {small} \[ \sum_{e,p} \!
    \underbrace {\vphantom {\Bigg|}
      \left[ \frac {3-r}{2}-(d-e)(r+1) \right]^{k(d)}_p \!\! (R^p)_a{}^b}
      _{\mbox {(A)}}
    \underbrace {\vphantom {\Bigg|}
      (-1)^{1-a-(d-e)(r+1)} \langle \tau_{\bullet}(T^a)
      \rangle_{0,d-e}}_{\mbox {(B)}}
    \underbrace {\vphantom {\Bigg|}
      \langle \tau_{\bullet}(T_b) T^{b_j} \rangle_{0,e-e_j}}
      _{\mbox {(C)}}. \] \end {small}%
  The three terms in this expression can all be expressed easily in terms of
  generating functions. Recalling that $ (R^p)_a{}^b = (r+1)^p \delta_{a+p,b}
  $, the (A) term is by definition equal to the $ z^p $-coefficient of
      \[ \delta_{a+p,b} \prod_{i=0}^{k(d)} \left(
           (r+1)z+\frac {3-r}2-(d-e)(r+1)+i \right). \]
  The (B) and (C) terms are rational 2-point invariants of $ \PP^r $ which have
  been computed in \cite {P} section 1.4: the Gromov-Witten invariant in (B)
  (without the sign) is equal to the $ z^a $-coefficient of $ \prod_{i=1}^{d-e}
  \frac 1 {(z+i)^{r+1}} $. So including the sign factor we get the $ z^a
  $-coefficient of $ - \prod_{i=1}^{d-e} \frac 1 {(z-i)^{r+1}} $. The (C) term
  is again by \cite {P} equal to the $ z^{r-b} $-coefficient of $
  \frac 1 {(z+e-e_j)^{b_j+1}} \prod_{i=1}^{e-e_j-1} \frac 1 {(z+i)^{r+1}} $.

  Multiplying these expressions and performing the sums over $a$, $b$, and $p$,
  we find that $ W_{d,j} $ is the $ z^r $-coefficient of
    \[ - \sum_e \frac {\prod_{i=0}^{k(d)} \left(
         (r+1)z+\frac {3-r}2-(d-e)(r+1)+i \right)}{
	 (z+e-e_j)^{b_j+1} \prod_{i=1}^{d-e} (z-i)^{r+1} \cdot
                           \prod_{i=1}^{e-e_j-1} (z+i)^{r+1}}, \]
  which can be rewritten as the sum of residues
    \[ - \sum_e \res_{z=0} \frac {\prod_{i=0}^{k(d)} \left(
         (r+1)z+\frac {3-r}2-(d-e)(r+1)+i \right)}{
	 (z+e-e_j)^{b_j+1} \prod_{i=e-d}^{e-e_j-1} (z+i)^{r+1}} \; dz. \]
  Note that this fraction depends on $z$ and $e$ only in the combination $ z+e
  $. Consequently, instead of summing the above residues at 0 over all $e$ we
  can as well set $ e=d $ and sum over all poles $ z \in \CC $ of the rational
  function. So we see that $ W_{d,j} $ is equal to
    \[ - \sum_{z_0 \in \CC} \res_{z=z_0} \frac {\prod_{i=0}^{k(d)} \left(
         (r+1)z+\frac {3-r}2+i \right)}{
	 (z+d-e_j)^{b_j+1} \prod_{i=0}^{d-e_j-1} (z+i)^{r+1}} \; dz. \]
  By the residue theorem this is nothing but the residue at infinity of
  our rational function. So we conclude that
    \[ W_{d,j} = \res_{z=0} \frac {\prod_{i=0}^{k(d)} \left(
         \frac {r+1}z+\frac {3-r}2+i \right)}{
	 (\frac 1z+d-e_j)^{b_j+1} \prod_{i=0}^{d-e_j-1} (\frac 1z+i)^{r+1}}
         \; d(\frac 1z). \]
  Finally note that by equations (\ref {dim-cond}) and (\ref {k-det}) we have
  the dimension condition
    \[ k(d)+1 = j+b_j+(d-e_j)(r+1), \]
  so multiplying our expression with $ z^{k(d)+1} $ in the numerator and
  denominator we get
    \[ W_{d,j} = - \res_{z=0} \frac {\prod_{i=0}^{k(d)} \left(
         r+1+\left( \frac {3-r}2+i\right) z \right)}{
	 z^{j+1} (1+(d-e_j)z)^{b_j+1} \prod_{i=0}^{d-e_j-1} (1+iz)^{r+1}}
         \; dz. \]
  This proves the lemma.
\end {proof}

To avoid unnecessary factors in the determinants, let us divide row $d$ of $W$
by the non-zero number $ - (r+1)^{k(d)+1} $ and call the resulting matrix $
\tilde W $. So we will now consider $ N \times N $ submatrices of $ \tilde
W=(\tilde W_{d,j}) $, obtained by picking the first $N$ columns of any $N$
rows, where $ \tilde W_{d,j} $ is the $ z^j $-coefficient of
\begin {equation} \label {w-def}
  \frac {\prod_{i=0}^{k(d)} \left(
    1+\left( \frac {3-r}{2r+2}+\frac {i}{r+1} \right) z \right)}{
    (1+(d-e_j)z)^{b_j+1} \prod_{i=0}^{d-e_j-1} (1+iz)^{r+1}}.
\end {equation}
The following technical lemma is the main step in computing their determinants.

\begin {lemma} \label {technical}
  Assume that we are given $ N,n \in \NN $, $ M \in \ZZ $, $ q,c \in \RR $,
  and distinct integers $ a_0,\dots,a_N $. Set
    \[ f(z) = \sum_{k=0}^N \left(
         \prod_{i \neq k} \frac 1 {a_k-a_i} \cdot
         \prod_{i=M}^{na_k} \left(
           1+\frac {c+i}{n} \, z
         \right) \cdot
         (1+a_k z)^q \cdot
         \prod_{i=a_k}^{-1} (1+iz)^n
       \right) \]
  as a formal power series in $z$.
  \begin {enumerate}
  \item For any $ i \ge 0 $ the $ z^i $-coefficient of $ f(z) $ is a
    symmetric polynomial in $ a_0,\dots,a_N $ of degree at most $ i-N $. (In
    particular, it is zero for $ i<N $.)
  \item The $ z^N $-coefficient of $ f(z) $ is equal to
    \[ \frac 1 {N!} \; \prod_{i=1}^N \left(
         c+q-N+\frac {n+1}2+i
       \right). \]
  \end {enumerate}
\end {lemma}

\begin {proof}
  In the following proof, we will slightly abuse notation and vary the
  arguments given explicitly for the function $f$. So if we e.g.\ want to study
  how $ f(z) $ changes if we vary $c$, we will write $ f(z) $ also as $ f(z,c)
  $, and denote by $ f(z,c+1) $ the function obtained from $ f(z) $ when
  substituting $c$ by $ c+1 $.

  (i): It is obvious by definition that $ f(z) $ is symmetric in the $ a_i $.
  We will prove the polynomiality and degree statements by induction on $N$.

  ``$ N=0 $'': In this case we have
    \[ f(z) = \prod_{i=M}^{na_0} \left( 1+\frac {c+i}{n} \, z \right) \cdot
         (1+a_0 z)^q \cdot \prod_{i=a_0}^{-1} (1+iz)^n. \]
  We have to show that the $ z^i $-coefficient of $ f(z) $ is a polynomial
  in $ a_0 $ of degree at most $i$. Note that this property is stable under
  taking products, so if we write $ f(z) = \prod_{j=0}^n f^{(j)}(z) $ with
  \begin {align*}
    f^{(0)}(z) &= \prod_{i=M}^{0} \left( 1+\frac {c+i}{n} \, z \right) \cdot
                  (1+a_0 z)^q \\
    \mbox {and} \quad
    f^{(j)}(z) &= \prod_{i=0}^{a_0-1} \left( 1 + \frac {c+j}{n} \,
                  \frac {z}{1+iz} \right) \quad \mbox {for $ 1 \le j \le n $}
  \end {align*}
  then it suffices to prove the statements for the $ f^{(j)} $ separately.
  But the statement is obvious for $ f^{(0)} $, so let us focus on $ f^{(j)} $
  for $ j>0 $. Note that
    \[ f^{(j)} (z,a_0+1) = f^{(j)} (z,a_0) \cdot \left(
       1 + \frac {c+j}{n} \, \frac {z}{1+a_0 z} \right). \]
  So if $ f_i $ denotes the $ z^i $-coefficient of $ f(z) $ we get
  \begin {equation} \label {fj}
    f^{(j)}_i (a_0+1) - f^{(j)}_i (a_0) =
         \frac {c+j}{n} \, \sum_{k=0}^{i-1} (-a_0)^k
           f^{(j)}_{i-1-k}(a_0).
  \end {equation}
  The statement now follows by induction on $i$: it is obvious that the
  constant $z$-term of $ f^{(j)}(z) $ is 1. For the induction step, assume that
  we know that $ f^{(j)}_i $ is polynomial of degree at most $i$ in $ a_0 $ for
  $ i=0,\dots,i_0-1 $. Then the right hand side of (\ref {fj}) is polynomial of
  degree at most $ i_0-1 $ in $ a_0 $, so $ f^{(j)}_{i_0} $ is polynomial of
  degree at most $ i_0 $. This completes the proof of the $ N=0 $ part of (i).

  ``$ N \to N+1 $'': Note that
  \begin {equation} \label {rec-rel}
    f(z,N+1,a_0,\dots,a_{N+1}) = \frac {
      f(z,N,a_0,\dots,a_N)-f(z,N,a_1,\dots,a_{N+1})}{a_0-a_{N+1}}.
  \end {equation}
  By symmetry we have $ f(z,N,a_0,\dots,a_N)=f(z,N,a_1,\dots,a_{N+1}) $ if $
  a_0=a_{N+1} $. Hence every $ z^i $-coefficient of this expression is a
  polynomial in the $ a_k $. Its degree is at most $ (i-N)-1 $ by the induction
  hypothesis. This proves (i).

  (ii): By (i) the $ z^N $-coefficient of $ f(z,N) $ does not depend on the
  choice of $ a_k $, so we can set $ a_k=a+k $ for all $k$ and keep only $a$ as
  a variable. It does not depend on $M$ either, as a shift $ M \mapsto M \pm 1
  $ corresponds to multiplication of $ f(z) $ with $ (1+\alpha z)^{\mp 1} $ for
  some $ \alpha $, which does not affect the leading coefficient of $ f(z) $.
  So we can set $ M=1 $ without loss of generality.

  The recursion relation (\ref {rec-rel}) now reads
  \begin {equation} \label {rec-rel-2}
    f(z,N+1,a) = \frac {f(z,N,a+1)-f(z,N,a)}{N+1}.
  \end {equation}
  By (i) the $ z^i $-coefficient of $ f(z) $ has degree at most $ i-N $ in $a$.
  So if we denote by $ f_i $ the $ a^{i-N} $-coefficient of the $ z^i
  $-coefficient of $ f(z,a) $, comparing the $ z^i $-coefficients in (\ref
  {rec-rel-2}) yields $ f_i(N+1) = \frac {i-N}{N+1} f_i(N) $ and therefore
    \[ f_N(N) = \frac {1}{n} \cdot \frac {2}{n-1} \cdots \frac {n}{1}
       \cdot f_N(0) = f_N(0). \]
  In other words, instead of computing the $ z^N $-coefficient of $ f(z,N) $ we
  can as well set $ N=0 $ and compute the $ a^N $-coefficient (i.e.\ the
  leading coefficient in $a$) of the $ z^N $-coefficient of $ f(z,N=0) $. So
  let us set $ N=0 $ to obtain
    \[ f(z) = \prod_{i=1}^{na} \left(
           1+\frac {c+i}{n} \, z
         \right) \cdot
         (1+a z)^q \cdot
         \prod_{i=a}^{-1} (1+iz)^n, \]
  and denote by $ g_N $ the $ a^N $-coefficient of the $ z^N $-coefficient of $
  f(z,a) $. Moreover, set $ g(z) = \sum_{N \ge 0} g_N z^N $. Our goal is then
  to compute $ g(z) $.

  We will do this by analyzing how $ f(z) $ (and thus $ g(z) $) varies when we
  vary $q$, $c$, or $n$. To start, it is obvious that
  \begin {equation} \label {change-q}
    g(z,q+\alpha) = (1+z)^\alpha g(z,q)
  \end {equation}
  for all $ \alpha \in \RR $. Next, note that
    \[ f(z,c+1) = f(z,c) \cdot \frac {1+\frac {c+1+na}{n}\,z}{
         1+\frac {c+1}{n}\,z}. \]
  For $ g(z) $ we can drop all terms in which the degree in $a$ is smaller than
  the degree in $z$. So we conclude
    \[ g(z,c+1) = (1+z) g(z,c). \]
  Combining this with (\ref {change-q}) we see that $ g(z) $ will depend on $q$
  and $c$ only through their sum $ q+c $. So in what follows we can set $ c=0
  $, and replace $q$ by $ q+c $ in the final result.

  Varying $n$ is more complicated. We have
    \[ f(z,n+1) = f(z,n) \cdot \underbrace {\prod_{j=0}^n \prod_{i=1}^a \left(
         1+\frac {j}{n(n+1)} \cdot \frac {z}{1+(i-\frac jn)z}
       \right)}_{=: \tilde f(z)}. \]
  Recall that for $ g(z) $ we only need the summands in $ \tilde f(z) $ in
  which the degree in $a$ is equal to the degree in $z$. So let us denote the $
  a^N $-coefficient of the $ z^N $-coefficient of $ \tilde f(z,a) $ by $ \tilde
  g_N $, and assemble the $ \tilde g_N $ into a generating function $ \tilde
  g(z) = \sum_{N \ge 0} \tilde g_N z^N $, so that $ g(z,n+1) = g(z,n) \cdot
  \tilde g(z) $. To determine $ \tilde g(z) $ compare the $ a^N $-coefficient
  of the $ z^N $-coefficient in the recursive equation
    \[ \frac {\tilde f(z,a) - \tilde f(z,a-1)}{z} =
       \tilde f(z,a-1) \cdot \frac 1z \, \left( \prod_{j=0}^n \left(
         1+\frac {j}{n(n+1)} \cdot \frac {z}{1+(a-\frac jn)z}
       \right) -1 \right). \]
  On the left hand side this coefficient is $ (N+1)\tilde g_{N+1} $. On the
  right hand side it is the $ z^N $-coefficient of
    \[ \tilde g(z) \cdot \left( \sum_{j=0}^N \frac j {n(n+1)} \right)
       \cdot \frac 1 {1+z} = \tilde g(z) \cdot \frac 12 \, \frac 1{1+z}. \]
  So we see that
    \[ \frac {d\tilde g(z)}{dz} = \frac 12 \, \frac 1{1+z} \, \tilde g(z). \]
  Together with the obvious initial condition $ \tilde g(0)=1 $ we conclude
  that $ \tilde g(z) = \sqrt {1+z} $, and therefore
    \[ g(z,n+1) = g(z,n) \cdot \sqrt {1+z}. \]
  Comparing this with (\ref {change-q}) we see that $ g(z) $ depends on $n$ and
  $q$ only through the sum $ q+\frac n2 $. We can therefore set $ n=1 $ and
  then replace $q$ by $ q+\frac {n-1}2 $ in the final result.

  But setting $n$ to 1 (and $c$ to 0) we are simply left with
    \[ f(z) = \prod_{i=1}^a (1+iz) \cdot
         (1+a z)^q \cdot
         \prod_{i=a}^{-1} (1+iz) = (1+az)^{q+1}. \]
  So it follows that $ g(z) = (1+z)^{q+1} $ and therefore
    \[ g_N = \binom {q+1}{N} = \frac 1 {N!} \; \prod_{i=1}^N (q+1-N+i). \]
  Setting back in the $c$ and $n$ dependence, i.e.\ replacing $q$ by $
  q+c+\frac {n-1}2 $, we get the desired result.
\end {proof}

We are now ready to compute our determinant.

\begin {proposition} \label {det-final}
  Let $ (\tilde W)_{d \ge \delta, j \ge 0} $ be the matrix defined in equation
  (\ref {w-def}). Pick $N$ distinct integers $ d_0,\dots,d_{N-1} $ with $ d_i
  \ge \delta $ for all $i$. Then the determinant of the $ N \times N $
  submatrix of $ \tilde W $ obtained by picking the first $N$ columns of rows $
  d_0,\dots,d_{N-1} $ is equal to
    \[ \frac {\prod_{i>j} (d_i-d_j)}{\prod_{i=1}^{N-1} i!} \cdot
       \prod_{i=1}^{N-1} \left( i+\frac 12 \right)^{N-i}. \]
  In particular, this determinant is never zero.
\end {proposition}

\begin {proof}
  We prove the statement by induction on $N$. The result is obvious for $ N=1 $
  as every entry in the first column (i.e.\ $ j=0 $) of $ \tilde W $ is equal
  to 1. So let us assume that we know the statement for a given value $N$. We
  will prove it for $ N+1 $.

  Denote by $ \Delta(d_0,\dots,d_{N-1}) $ the determinant of the $ N \times N $
  submatrix of $ \tilde W $ obtained by picking the first $N$ columns of rows
  $ d_0,\dots,d_{N-1} $. Then by expansion along the last column and the
  induction assumption we get
  \begin {align*}
    \Delta(d_0,\dots,d_N) &=
      \sum_{k=0}^N (-1)^{k+N} \tilde W_{d_k,N} \cdot
        \Delta(d_0,\dots,d_{k-1},d_{k+1},\dots,d_N) \\
    &= \frac {\prod_{i>j} (d_i-d_j)}{\prod_{i=1}^{N-1} i!} \cdot
        \prod_{i=1}^{N-1} \left( i+\frac 12 \right)^{N-i} \cdot \sum_{k=0}^N
	\left( \prod_{i \neq k} \frac {1}{d_k-d_i}
          \cdot \tilde W_{d_k,N} \right).
  \end {align*}
  But by lemma \ref {technical}, applied to the values $ n=r+1 $, $
  a_k=d_k-e_N $, $ q=-b_N-1 $, $ M=(r+1)(d_k-e_N)-k(d_k) = 1-N-b_N $, and
  $ c=\frac {3-r}{2}-M = \frac {1-r}{2} +N+b_N $, the sum in this expression is
  equal to
    \[ \frac 1{N!} \; \prod_{i=1}^N \left( i+\frac 12 \right). \]
  Inserting this into the expression for the determinant, we obtain
    \[ \Delta(d_0,\dots,d_N) =
         \frac {\prod_{i>j} (d_i-d_j)}{\prod_{i=1}^{N} i!} \cdot
       \prod_{i=1}^{N} \left( i+\frac 12 \right)^{N+1-i}, \]
  as desired.
\end {proof}

\begin {remark}
  It should be remarked that the expression for the determinant in proposition
  \ref {det-final} is surprisingly simple, given the complicated structure of
  the Virasoro conditions and the topological recursion relations. It would be
  interesting to see if there is a deeper relation between these two sets of
  equations that is not yet understood and explains the simplicity of our
  results.
\end {remark}

Combining the arguments of this section, we see that the systems of linear
equations obtained from the Virasoro conditions and the topological recursion
relations in section \ref {virasoro} are always solvable. This completes the
proof of theorem \ref {reconstruction}.


\section {Some numbers} \label {numbers}

In this section we will give some examples of invariants that have been
computed using the method of this paper.

\begin {example} (The Caporaso-Harris numbers, see \cite {CH})
  The following table shows some numbers of curves in $ \PP^2 $ of genus $g$
  and degree $d$ through $ 3d-1+g $ points, i.e.\ the Gromov-Witten invariants
  $ \langle \pt^{3d-1+g} \rangle_{g,d} $. They can be found either by the
  Caporaso-Harris method or by applying the techniques of this paper.
    \[ \begin {array}{|c|rrrrrrr|} \hline
        & d=1 & d=2 & d=3 & d=4 &   d=5 &      d=6 &          d=7 \\ \hline
    g=0 &   1 &   1 &  12 & 620 & 87304 & 26312976 &  14616808192 \\
    g=1 &   0 &   0 &   1 & 225 & 87192 & 57435240 &  60478511040 \\
    g=2 &   0 &   0 &   0 &  27 & 36855 & 58444767 & 122824720116 \\
    g=3 &   0 &   0 &   0 &   1 &  7915 & 34435125 & 153796445095 \\
    g=4 &   0 &   0 &   0 &   0 &   882 & 12587820 & 128618514477 \\ \hline
  \end {array} \]
\end {example}

\begin {example}
  We list some 1-point invariants of $ \PP^2 $, i.e.\ invariants of the form
  $ \langle \tau_m(\gamma) \rangle_{g,d} $, where $m$ is determined by the
  dimension condition $ 3d+g=m+\deg \gamma $.
  \begin {small} \[ \begin {array}{|@{\;}c@{\;}|@{\;}r@{\;}r@{\;}r@{\;}|
                                             @{\;}r@{\;}r@{\;}r@{\;}|
                                             @{\;}r@{\;}r@{\;}r@{\;}|} \hline
         & & \gamma=\pt & & & \gamma=H & & & \gamma=1 & \\
	 & d=0 & d=1 & d=2 & d=0 & d=1 & d=2 & d=0 & d=1 & d=2 \\ \hline
      g=0 & - & 1 & \frac {1}{8}
          & - & -3 & -\frac {9}{16}
          & - & 6 & \frac {3}{2} \\
      g=1 & - & 0 & \frac {1}{32}
          & -\frac {1}{8} & \frac {1}{8} & -\frac {3}{32}
          & \frac {1}{8} & -\frac {1}{4} & \frac {23}{128} \\
      g=2 & 0 & -\frac {1}{240} & -\frac {1}{960}
          & -\frac {1}{960} & -\frac {1}{960} & \frac {13}{1536}
          & \frac {7}{640} & \frac {1}{128} & -\frac {27}{1280} \\
      g=3 & 0 & \frac {1}{3360} & -\frac {1}{16128}
          & -\frac {1}{40320} & 0 & -\frac {163}{645120}
          & \frac {41}{161280} & -\frac {97}{161280} & \frac {43}{36864} \\
      g=4 & 0 & -\frac {1}{80640} & \frac {11}{1075200}
          & -\frac {1}{1075200} & -\frac {1}{153600} & -\frac {1}{147456}
          & \frac {127}{12902400} & \frac {173}{4300800} & -\frac
	    {4567}{103219200} \\ \hline
       \end {array} \] \end {small}%
\end {example}

\begin {example}
  The following table gives some Gromov-Witten invariants $ \langle \pt^{2d}
  \rangle_{g,d} $ of $ \PP^3 $, i.e.\ the virtual number of genus-$g$
  degree-$d$ curves in $ \PP^3 $ through $ 2d $ points. Note that these
  Gromov-Witten invariants are not enumerative for $ g>0 $. The precise
  relationship between the Gromov-Witten invariants and the enumerative
  numbers is not yet known.
  \[ \begin {array}{|c|rrrrrrr|} \hline
        & d=1 & d=2 & d=3 & d=4 & d=5 &  d=6 &    d=7 \\ \hline
    g=0 &   1 &   0 &   1 &   4 & 105 & 2576 & 122129 \\
    g=1 & -\frac {1}{12}
              &   0
                    & -\frac {5}{12}
                          & -\frac {4}{3}
                                & -\frac {147}{4}
                                      & \frac {1496}{3}
                                             & \frac {1121131}{12} \\
    g=2 & \frac {1}{360}
              & 0
                    & \frac {1}{12}
                          & -\frac {1}{180}
                                & -\frac {49}{8}
                                      & -\frac {7427}{5}
                                             & -\frac {4905131}{45} \\
    g=3 & -\frac {1}{20160}
              & 0
                    & -\frac {43}{4032}
                          & \frac {103}{1080}
                                & \frac {473}{64}
                                      & \frac {206873}{270}
                                             & \frac {283305113}{8640} \\
    g=4 & \frac {1}{1814400}
              & 0
                    & \frac {713}{725760}
                          & -\frac {26813}{907200}
                                & -\frac {833}{320}
                                      & -\frac {12355247}{56700}
                                             & -\frac {1332337}{34560}\\ \hline
  \end {array} \]
\end {example}


\end {document}